\documentclass[11pt]{article}
\usepackage{amsmath,amsfonts,epsfig, psfrag}
\newcommand{\bb}{\mathbb}

\newcommand{\cx}{{\bb C}}

\newcommand{\integers}{{\bb Z}}

\newcommand{\reals}{{\bb R}}

\newcommand{\hthree}{{\bb H}^3}

\newcommand{\widemargins}{
\setlength{\textwidth}{6.0in}
\setlength{\oddsidemargin}{0.25in}
\setlength{\evensidemargin}{0.25in}
}

\newcommand{\qed}[1]{\nopagebreak[4]{\tiny \hfill
\fbox{\ref{#1}} \linebreak }\pagebreak[2]}
\newcommand{\del}{\partial}

\newcommand{\chat}{\widehat{\cx}}

\newcommand{\id}{\operatorname{id}}

\newtheorem{theorem}{Theorem}[section]

\newtheorem{lemma}[theorem]{Lemma}
\newtheorem{cor}[theorem]{Corollary}
\newtheorem{conj}[theorem]{Conjecture}

\newcommand{\cC}{{\cal C}}

\widemargins

\psfrag{S}{$\Sigma_K$}
\psfrag{S1}{$\Sigma^+_K$}
\psfrag{S2}{$\Sigma^-_K$}
\psfrag{c1}{$c_1$}
\psfrag{c2}{$c_2$}
\psfrag{c}{$c$}
\psfrag{X0}{$X_0$}
\psfrag{X1}{$X_1$}
\psfrag{X2}{$X_2$}
\psfrag{A1}{$A_1$}
\psfrag{A2}{$A_2$}
\psfrag{B1}{$B_1$}
\psfrag{B2}{$B_2$}
\psfrag{ap}{$A^+_1$}
\psfrag{a-}{$A^-_1$}
\psfrag{Ap}{$A^+_2$}
\psfrag{A-}{$A^-_2$}
\psfrag{bp}{$B^+_1$}
\psfrag{b-}{$B^-_1$}
\psfrag{Bp}{$B^+_2$}
\psfrag{B-}{$B^-_2$}
\psfrag{MZ}{$M_\integers - A_\integers$}
\psfrag{M}{$M - A$}
\psfrag{MC}{$M_c$}
\psfrag{H}{$H-(B_1 \cup B_2)$}

\begin{document}
\title{Projective structures with degenerate holonomy and the Bers
density conjecture}

\author{K. Bromberg\footnote{This work was partially supported by a
grants from the NSF and the Clay Mathematics Institute}}

\date{\today}

\maketitle

\begin{abstract}
\noindent
We prove the Bers' density conjecture for singly degenerate Kleinian surface groups without parabolics.
\end{abstract}

\section{Introduction}
In this paper we address a conjecture of Bers about singly degenerate
Kleinian groups. These are discrete subgroups of $PSL_2\cx$ that
exhibit some unusual behavior:
\begin{itemize}
\item As groups of projective transformations of the
Riemann sphere $\chat$ they act properly discontinuously on a
topological disk whose closure is all of $\chat$.
\item As groups of hyperbolic isometries their action on $\hthree$ is not convex co-compact.
\item Viewed as a dynamical systems they are not structurally stable. 
\end{itemize}
These groups were first discovered by Bers in his paper \cite{Bers:bdry} where he made the conjecture that will be the
focus of our work here.

Let $M= S \times [-1,1]$ be an $I$-bundle over a closed surface $S$ of
genus $>1$. We will be interested in the space $AH(S)$ of all Kleinian groups isomorphic to
$\pi_1(S)$. By a theorem of Bonahon, this is equivalent to studying
complete hyperbolic structures on the interior of $M$. A generic hyperbolic structure
on $M$ is {\em quasi-fuchsian} and the geometry is well understood
outside of a compact set.
In particular, although the geometry of the surfaces $S \times \{t\}$
will grow exponentially as $t$ limits to $-1$ or $1$, the conformal
structures will stabilize and limit to Riemann surfaces $X$ and
$Y$. Then $M$ can be conformally compactified by viewing $X$ and
$Y$ as conformal structures on $S \times \{-1\}$ and $S \times \{1\}$,
respectively. Bers showed that $X$ and $Y$ parameterize the space
$QF(S)$ of all quasi-fuchsian structures. In other words $QF(S)$ is 
isomorphic to $T(S) \times T(S)$ where $T(S)$ is the Teichm\"uller
space
of marked conformal structures on $S$. Let the {\em Bers slice} $B_X$ be the slice of
$QF(S)$ obtained by fixing $X$ and letting $Y$ vary in $T(S)$. 

This gives an interesting model of $T(S)$ because $B_X$ naturally
embeds as a bounded domain in the space $P(X)$ of projective
structures on $S$ with conformal structure $X$. The closure
$\overline{B}_X$ of $B_X$ in $P(X)$ is then a compactification of
Teichm\"uller space. A point in $\del B_X = \overline{B}_X - B_X$ will
again correspond to a complete hyperbolic structure on $M$. As
with structures in $B_X$, the surfaces $S \times \{t\}$ will converge
to the conformal structure $X$ as $t \rightarrow -1$. However, as $t
\rightarrow 1$ the structures will not converge. 

There are three possibilities for the limiting geometry of the $S
\times 
\{t\}$. In the simplest case there will be an essential simple closed
curve (or a collection of curves) $c$ on $S$ such that the length of
$c$ on $S \times \{t\}$ limits to zero, while on the complement of $c$
the surfaces grow exponentially but converge to a cusped conformal
structure. In this case $M$ is {\em geometrically finite}. In the
other case there will be a sequence $t_i \rightarrow 1$ such that $S
\times \{t_i\}$ has bounded area yet for any simple closed curve
$c$ on $S$ the length of $c$ on $S \times \{t_i\}$ will go to
infinity as $t_i \rightarrow 1$. In other words, the geometry of the
$S \times \{t_i\}$ is bounded but still changing radically. Such
manifolds are {\em singly degenerate}. The final possibility is
that $M$ may have a combination of the first two behaviors. 

Understanding such structures is a motivating problem in hyperbolic
3-manifolds and Kleinian groups. Bers made the following conjecture: 

\begin{conj}[Bers Density Conjecture \cite{Bers:bdry}]
\label{bersconj}
Let $\Gamma \in AH(S)$ be a Kleinian group. If $M = \hthree/\Gamma$ is
singly degenerate then $\Gamma \in \overline{B}_X$ where $X$ is the
conformal boundary of $M$.
\end{conj}

There are some special cases where the conjecture is known. Abikoff
\cite{Abikoff:degenerating} proved the conjecture when $M$ is
geometrically finite. Recently Minsky \cite{Minsky:bounded} has
proved the conjecture in the case where there is a lower bound on the
length of any closed geodesic in $M$ and $\Gamma$ has no parabolics 
($M$ has {\em bounded geometry}). In a separate, earlier paper
(\cite{Minsky:torus}), Minsky also proved the conjecture if $S$ is a
punctured torus. In this paper we prove the conjecture when $M$ has a
sequence of closed geodesics $c_i$ whose length limits to zero ($M$
has {\em unbounded geometry}).  Combined with Minsky's result we have
an almost complete resolution of Bers's conjecture:

\medskip
\noindent
{\bf Theorem \ref{bers}}
{\em Assume that $\Gamma \in AH(S)$ has no parabolics.
If $M = \hthree/\Gamma$ is singly degenerate then $\Gamma \in
\overline{B}_X$ where $X$ is the conformal boundary of $M$.}

\medskip
\noindent
There is a more general version of the density conjecture due to Sullivan and Thurston. It states that every finitely generated Kleinian group is an algebraic limit of geometrically finite Kleinian groups. In joint work with Brock (\cite{Brock:Bromberg:density}) we use some of the ideas of this paper to prove this more general conjecture for freely indecomposable Kleinian groups without parabolics.

The condition
that $\Gamma$ has no parabolics is a technical one and we believe with
more work present techniques could be used to prove the complete
conjecture. More precisely, if the surface $S$ has punctures then
instead of studying all Kleinian groups isomorphic to $\pi_1(S)$,
$AH(S)$ is the space of Kleinian groups in which all of the punctures
are parabolic. If one could prove Conjecture \ref{bersconj} for all
$\Gamma \in AH(S)$ such that all parabolics in $\Gamma$ correspond to
punctures then the entire conjecture would follow. If $M = \hthree/\Gamma$
has unbounded geometry most of the work in this paper generalizes
easily. If $M$ has bounded geometry then one needs to generalize
Minsky's work. In particular, most of Minsky's work applies in this
setting; it is only his earliest paper on the problem
(\cite{Minsky:teichmuller}) that needs to be generalized.

We also remark that the density conjecture is a consequence of the
ending lamination conjecture. In fact, Minsky's results on the density
conjecture are a consequence of his work on the ending lamination
conjecture. More recently Brock, Canary and Minsky have annouced work that completes Minsky program to prove the full ending lamination conjecture (\cite{Minsky:lipschitz, Brock:Canary:Minsky:elc}).

We now outline our results.

Our approach to Conjecture \ref{bersconj} is to understand projective
structures with singly
degenerate holonomy. Our study will be guided by Goldman's
classification of all projective structures with quasi-fuchsian
holonomy. In particular, the two conformal structures $X$ and $Y$ that
compactify a quasi-fuchsian manifold also have projective structures
$\Sigma^-$ and $\Sigma^+$. Goldman showed that all projective
structures with quasi-fuchsian holonomy are obtained by {\em grafting}
on $\Sigma^-$ or $\Sigma^+$. For a singly degenerate group we still
have the projective structure $\Sigma^-$ and all of its graftings. On
the other hand, while the projective structure $\Sigma^+$ is gone we
will show that its graftings still exist.

We will use these projective structures to construct a family of
quasi-fuchsian hyperbolic cone-manifolds that converge to the singly
degenerate manifold $M$. Here is our main construction. By a theorem of Otal
\cite{Otal:short}, any sufficiently short geodesic $c$ will be 
unknotted. That is the product structure can be chosen such that $c$
is a simple closed curve on $S \times \{0\}$. Let $A$ be the annulus
$c \times [0,1)$ and let $A_\integers$ be a lift of $A$ to the
$\integers$-cover $M_\integers$ of $M$ associated to $c$. Now remove
$A$ from $M$ and $A_\integers$ from $M_\integers$ and take the metric
completion of both spaces. Both of these spaces will be manifolds with
boundary isometric to two copies of $A$ meeting at the geodesic
$c$. Next, glue the two manifolds together along their isometric
boundary to form a new manifold $M_c$. This new manifold will be
homeomorphic to $M$ but the hyperbolic structure will be singular
along the geodesic $c$. In particular $M_c$ will be a hyperbolic {\em
  cone-manifold}, for a cross section of a tubular neighborhood of $c$
will be a cone of cone angle $4\pi$. We will show:

\medskip
\noindent
{\bf Theorem \ref{coneman}} {\em
The hyperbolic cone-manifold $M_c$ is a quasi-fuchsian cone-manifold
with projective boundary $\Sigma$ and $\Sigma_c$.}

\medskip
The lower half of $M_c$ is isometric to the lower half of $M$ and is
therefore compactified by the same projective structure $\Sigma$ on the
conformal structure $X$. The upper half of $M_c$ will be compactified
by the new projective structure $\Sigma_c$ which will have conformal
structure $Y_c$. Then there is a unique quasi-fuchsian group
$\Gamma_c \in B_X$ such that $M'_c = \hthree/\Gamma_c$ has conformal
boundary $X$ and $Y_c$.

If $M$ has unbounded geometry there will be a sequence of closed
geodesics $c_i$ with $\operatorname{length}(c_i) \rightarrow
0$. Repeating the above construction for each $c_i$ we obtain cone
manifolds $M_i$ and quasi-fuchsian manifolds $M'_i =
\hthree/\Gamma_i$. Let $\Sigma_i$ be the component of the projective
boundary of $M'_i$ corresponding to $X$. The final step is to bound the distance between $\Sigma_i$ and $\Sigma$ in terms of $\operatorname{length}(c_i)$.

This is done using the deformation theory
of hyperbolic cone-manifolds developed by Hodgson and Kerckhoff for closed manifolds and extended by the author to geometrically finite cone-manifolds. For
each $M_i$ we can use this deformation theory to find a smooth one parameter family of cone manifolds
that interpolates between $M_i$ and $M'_i$. Furthermore this
deformation theory allows us to control how the projective structure
$\Sigma$ deforms to the projective structure $\Sigma_i$. In
particular there is a canonical way to define a metric on $P(X)$ and
in this metric we have:
$$d(\Sigma, \Sigma_i) \leq K\operatorname{length}(c_i).$$
Therefore $\Sigma_i \rightarrow
\Sigma$ in $P(X)$ which implies that $\Gamma_i \rightarrow \Gamma$ in
$AH(S)$ and $\Gamma \in \overline{B}_X$. A novel feature of the above estimate is its use of the analytic theory of cone-manifolds to obtain results about infinite volume hyperbolic 3-manifolds. This approach has turned out to be fruitful in other problems (see \cite{Bromberg:schwarz, Brock:Bromberg:density}, \cite{Brock:Bromberg:Evans:Souto}) and we expect it will have further applications as well.

\begin{figure}[h]
\begin{center}
\epsfig{file=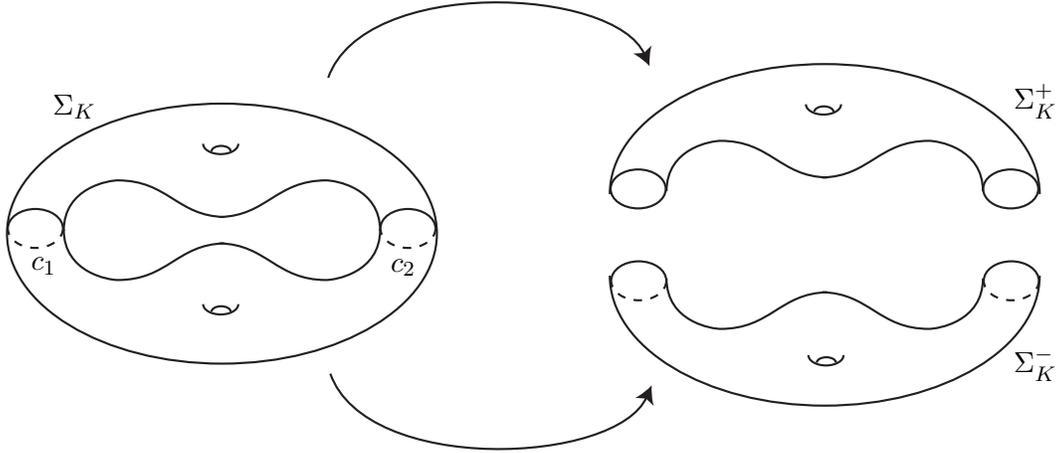, scale=.8}
\caption{Cutting $\Sigma_K$ along $c_1$ and $c_2$ produces $\Sigma^+_K$ and $\Sigma^-_K$.}
\label{handlebody}
\end{center}
\end{figure}

\begin{figure}[h]
\begin{center}
\epsfig{file=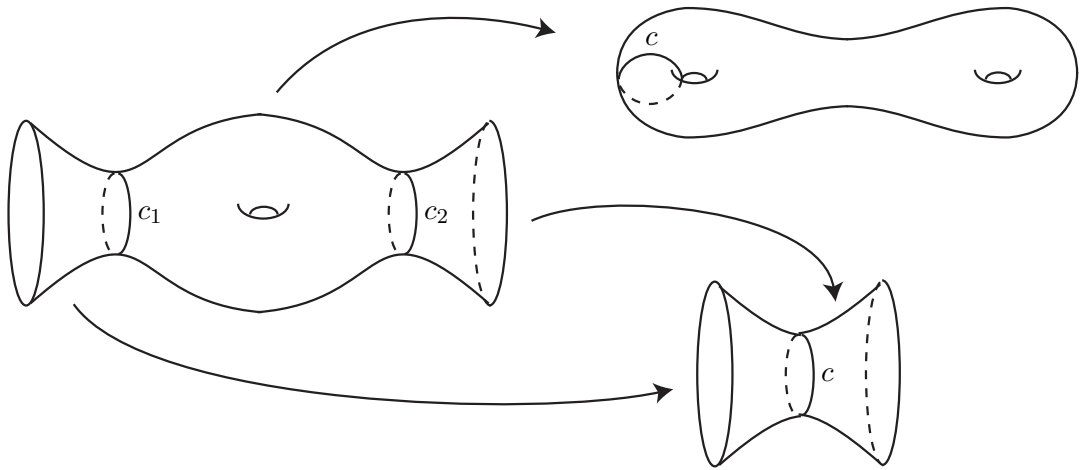, scale=.85}
\caption{If we cut $S'$ along $c_1$ and $c_2$ we have three pieces which can be reglued to form the original surface $S$ and the cover of $S$ associated to $c$.}
\label{cover}
\end{center}
\end{figure}

{\bf Acknowledgments.} The author would like to thank Manny Gabet for drawing Figures \ref{handlebody} and \ref{cover} and Jeff Brock for many helpful comments on a draft version of this paper.

\section{Preliminaries}
\subsection{Kleinian groups}
A Kleinian group $\Gamma$ is a discrete subgroup of $PSL_2\cx$. In
this paper we will assume that all Kleinian groups are torsion
free. The Lie group $PSL_2\cx$
acts as both projective transformations of the Riemann sphere
$\chat$ and as isometries on hyperbolic 3-space $\hthree$. The
union $\hthree \cup \chat$ is naturally topologized as a closed
3-ball such that the action of $PSL_2\cx$ on $\hthree$ extends
continuously to the action on $\chat$.

The domain of discontinuity $\Omega \subset \chat$ for $\Gamma$ is the largest subset of
$\chat$ such that $\Gamma$ acts properly discontinuously. The limit
set $\Lambda = \chat -\Omega$ is the complement of $\Omega$ in
$\chat$. The group $\Gamma$ will act properly discontinuously on all
of $\hthree$ so the quotient $\hthree/\Gamma$ will be a
3-manifold. The quotient $(\hthree \cup \Omega)/\Gamma$ will be a
3-manifold with boundary.

\subsection{Projective structures} Let $S$ be a surface. A {\em
projective structure} $\Sigma$ on $S$ is an atlas of charts to
$\chat$ with transition maps elements of $PSL_2\cx$, the group of projective transformations of $\chat$. If $\Gamma$
is a Kleinian group isomorphic to $\pi_1(S)$ and $\Omega'$ is a connected
component of $\Omega$ that is fixed by $\Gamma$ then the quotient
$\Omega'/\Gamma$ will be a projective structure on $S$.

As projective transformations are conformal maps, a projective
structure $\Sigma$ also defines a conformal structure $X$ on $S$. If
$T(S)$ is the Teichm\"uller space of marked conformal structures on
$S$ 
and $P(S)$ is the space of projective structures, then there is a map
$P(S) \longrightarrow T(S)$ defined by $\Sigma \mapsto X$.

Let $P(X)$ be the pre-image of $X$ in $P(S)$ under this map. It is
well known that $P(X)$ can be identified with the vector space of holomorphic
quadratic differentials on $X$. In particular, two projective
structures $\Sigma, \Sigma' \in P(X)$ differ by a holomorphic
quadratic differential $\phi$. If $\rho$ is the hyperbolic metric on $X$ then $\phi \rho^{-2}$ is a function on $X$. We let $\|\phi\|_\infty$ be the sup norm of this function. This defines a metric on $P(X)$ by 
$$d(\Sigma, \Sigma') = \|\phi \|_\infty.$$

A projective structure is {\em Fuchsian} if its is the quotient of a round disk in $\chat$. There is a unique Fuchsian element $\Sigma_F$ in $P(X)$ and we let $\|\Sigma\|_\infty = d(\Sigma, \Sigma_F)$.

\subsection{Hyperbolic structures} 
A {\em hyperbolic structure} on a 3-manifold $M$ is a Riemannian
metric with constant sectional curvature equal to $-1$. Equivalently,
a hyperbolic structure can be defined as an atlas of local charts to
$\hthree$ with transition maps hyperbolic isometries.

We will also be interested in certain singular hyperbolic
structures. We let $\hthree_\alpha$ be $\reals^3$ with cylindrical
coordinates $(r, \theta, z)$ and the Riemannian metric $dr^2 +
\sinh^2r d\theta^2 + \cosh^2 r dz^2$ where $\theta$ is measure modulo
$\alpha$. The metric on $\hthree_\alpha$ is a smooth metric of
constant sectional curvature $\equiv -1$ when $r \neq 0$. It extends
to a complete, singular metric on all of $\hthree_\alpha$. The
sub-surfaces where $z$ is constant are hyperbolic
planes away from $r=0$. At $r=0$ there is a cone-singularity with cone
angle $\alpha$.

If $\alpha = 2\pi$ then $\hthree_\alpha$ is isometric to $\hthree$. If
$\alpha={2\pi n}$ where $n$ is a positive integer then there is an
obvious map from $\hthree_\alpha$ to $\hthree$ that is a local
isometry when $r \neq 0$ and has an order $n$ branch locus at $r=0$.

A metric on $M$ is a hyperbolic {\em cone-metric} if all points in $M$
are either modeled on $\hthree$ or the point $(0,0,0)$ in
$\hthree_\alpha$ for some $\alpha$. All points of the second type are
the singular locus $\cC$ for $M$. Clearly $\cC$ will consist of a
collection of disjoint, simple curves and all points in a component
$c$ of $\cC$ will be modeled on $\hthree_\alpha$ for some fixed
$\alpha$. Then $\alpha$ is the {\em cone-angle} for $c$. In this paper
we will assume that the singular locus consists of a finite collection
of simple closed curves.

\subsection{Kleinian surface groups}
The space of representations of $\pi_1(S)$ in $PSL_2\cx$ has a natural topology given be convergence on generators.
Let $AH(S)$ be the space of conjugacy classes of discrete, faithfull representations of $\pi_1(S)$ in $PSL_2\cx$ with the quotient topology. The image of each representation is a marked Kleinian group so we can view $AH(S)$ as a space of Klienian groups. A group $\Gamma \in AH(S)$ is
{\em quasi-fuchsian} if the limit set of $\Gamma$ is a Jordan curve.
The domain of discontinuity is then two topological disks $\Omega^-$
and $\Omega^+$. Let $X = \Omega^-/\Gamma$ and $Y = \Omega^+/\Gamma$ be
the quotient conformal structures on $S$. The assignment
$$\Gamma \mapsto (X,Y)$$
defines a map from the space of quasi-fuchsian structures $QF(S)$ to $T(S) \times T(S)$.

\begin{theorem}[Bers]
\label{qfparam}
The space $QF(S)$ of quasi-fuchsian structures is parameterized by
$T(S) \times T(S)$.
\end{theorem}

We define a {\em Bers' slice} by $B_X = \{X\} \times T(S) \subset
QF(S)$. This set of quasi-fuchsian groups is isomorphic to
$T(S)$. Bers 
showed that $B_X$ embeds as a bounded domain in $P(X)$ and therefore
the closure $\overline{B}_X$ is a compactification of Teichm\"uller
space. 

To understand a general $\Gamma \in AH(S)$, we need the following
important theorem:

\begin{theorem}[Bonahon \cite{Bonahon:tame}]
\label{bonahon}
The quotient 3-manifold $\hthree/\Gamma$ is homeomorphic to $S
\times (-1,1)$.
\end{theorem}

Bers original study was of groups $\Gamma \in AH(S)$ such that the
hyperbolic structure $\hthree/\Gamma$ on $S \times (-1,1)$ extends to
a projective structure $\Sigma$ on $S \times \{-1\}$.
If such a $\Gamma$ is not quasi-fuchsian and has no parabolics then $\Gamma$
is {\em singly degenerate}. For a singly degenerate group the domain
of discontinuity will be a single topological disk. On the other hand, if $\Gamma$ has
parabolics then they will correspond to a collection of disjoint,
essential, simple closed curves on $S$. The subgroups of $\Gamma$
corresponding to the components of the complement of the simple closed curves
will either be quasi-fuchsian groups or singly degenerate groups. We
will not investigate groups with parabolics in this paper.

There is a further dichotomy for hyperbolic 3-manifolds with
degenerate ends. Namely, $M$ has {\em bounded geometry} if there is a
lower bound on the length of any closed geodesic in $M$. Otherwise $M$
has {\em unbounded geometry}. As mentioned 
in the introduction, Minsky has proven Bers' conjecture (Conjecture
\ref{bersconj}) if $M$ has bounded geometry. In fact he has proven a much stronger result which we
only partially state here: 

\begin{theorem}[Minsky \cite{Minsky:bounded}]
\label{minsky}
Suppose $\Gamma \in AH(S)$ has no parabolics. Then if $M =
\hthree/\Gamma$ has bounded geometry, $\Gamma \in
\overline{QF(S)}$. Furthermore if $M$ is singly degenerate with
conformal boundary $X$ then $\Gamma \in \overline{B}_X$. 
\end{theorem}

\subsection{Quasi-fuchsian cone-manifolds}
There is an alternate definition of a quasi-fuchsian manifold that
extends naturally to cone-manifolds. A hyperbolic structure on the
interior of $S \times [-1,1]$ is quasi-fuchsian if it extends to a
projective structure on $S \times \{-1\}$ and $S \times \{1\}$. More
explicitly, for each point $x$ in $S \times \{-1\}$ or $S \times
\{1\}$ there exists a local chart from a neighborhood of $x$ in $S
\times [-1,1]$ (not simply a neighborhood in $S \times \{\pm1\}$) to
$\hthree \cup \chat$. The transition maps will again be elements of
$PSL_2\cx$ which act as automorphisms of $\hthree \cup \chat$. This definition agrees with our previous definition
of a quasi-fuchsian structure extends to a definition quasi-fuchsian hyperbolic cone-metrics on $S \times (-1,1)$.

\subsection{Handlebodies and Schottky groups}
A Kleinian group $\Gamma$ is a {\em Schottky group} if $H=(\hthree \cup \Omega)/\Gamma$ is a closed handlebody with boundary. A handlebody has many distinct product structures. In particular if $Y$ is a properly embedded surface in $H$ such that the inclusion map is a homotopy equivalence then $H$ is homeomorphic to $S \times [-1,1]$ with $S \times \{0\} = Y$.

\subsection{Grafting}
A projective structure $\Sigma$ on $S$ defines a {\em holonomy
  representation} of $\pi_1(S)$ via a {\em developing map}. In
particular, $\Sigma$ lifts to a projective structure $\tilde{\Sigma}$
on $\tilde{S}$. Any chart for $\Sigma$ will lift to a chart for
$\tilde{\Sigma}$. Since $\tilde{\Sigma}$ is simply connected, this
chart will extend to a projective map $D: \tilde{S} \longrightarrow \chat$ on all of
$\tilde{\Sigma}$. Furthermore there will be a representation $\rho:
\pi_1(S) \longrightarrow PSL_2\cx$ such that 
$$D(g(x)) = \rho(g)D(x)$$
for all $g \in \pi_1(S)$.
Then $D$ is a developing map with holonomy $\rho$. Note that $D$ is
unique up to post-composition with elements of $PSL_2\cx$ while $\rho$
is unique up to conjugacy. 

Now let $c$ be an essential, simple closed curve on $S$ and
$\tilde{c}$ a component of the pre-image of $c$ in $\tilde{S}$. Let $g \in \pi_1(S)$
generate the $\integers$-subgroup that fixes $\tilde{c}$. We also
assume that $\rho(g)$ is hyperbolic and that $D(\tilde{c})$ is a
simple arc in $\chat$. Then the quotient of $\chat$ minus the fixed
points of  
$\rho(g)$ is a torus $T$, $D(\tilde{c})$ descends to an essential
simple closed curve $c'$ on $T$ and $A = T - c'$ is a projective
structure on an annulus. We can form a new projective structure on $S$ by
removing the curve $c$ from the projective structure $\Sigma$ and glueing
in $n$ copies of $A$. The new projective structure is then a {\em
grafting} of $\Sigma$ along the curve $c$. Most importantly for our purposes the grafted projective structure has the same holonomy as
$\Sigma$.

Goldman used grafting to classify projective structures with
quasi-fuchsian holonomy. Let $\Gamma$ be a quasi-fuchsian group with
$\Omega^-$ and $\Omega^+$ the two components of the domain of
discontinuity. Then $\Sigma^{\pm} = \Omega^{\pm}/\Gamma$ are
projective structures on $S$. 

\begin{theorem}[Goldman \cite{Goldman:graft}]
\label{graft}
All projective structures with holonomy $\Gamma$ are obtained by
grafting on either $\Sigma^-$ or $\Sigma^+$.
\end{theorem} 

In the next section, we will conjecture that a similar classification holds for singly degenerate Kleinian groups.

\section{Projective structures}

Let $S$ be a closed surface of genus $g > 1$ and $\Gamma$ a singly
degenerate Kleinian group isomorphic to $\pi_1(S)$. Let $\Sigma =
\Omega/\Gamma$ be the quotient projective structure on $S$.

Let $D: \tilde{S} \longrightarrow \Omega \subset \chat$ be a
developing 
map for $\Sigma$ with holonomy representation $\rho: \pi_1(S)
\longrightarrow \Gamma$. Choose an essential simple closed curve $c$
on $S$ and let $\tilde{c}$ be the pre-image of $c$ in the universal
cover $\tilde{S}$. We will begin by assuming that $c$ is
non-separating and deal with the general case at the end of the
section. We also choose a component $\tilde{K}$ of $\tilde{S} -
\tilde{c}$. Note that since $c$ is non-separating the action of $\pi_1(S)$ on the components of $\tilde{S} - \tilde{c}$ has a single orbit. Let $\tilde{c}_K$ be
the components of $\tilde{c}$ which lie on the boundary of $\tilde{K}$.

Let $\Gamma_K$ be the subgroup of $\Gamma$ which fixes $D(\tilde{K})$. Then
$\Omega/\Gamma_K$ will be a cover of $\Sigma$ corresponding to the
restriction of $\pi_1(S)$ to $S - c$. In particular $\Gamma_K$ will be
isomorphic to $\pi_1(S-c)$, a free group on $2g-1$ generators. We also
note that $D(\tilde{c}_K)$ will descend to two simple closed curves
$c_1$ and $c_2$ on the cover $\Omega/\Gamma_K$.

Let $\Omega_K$ be the domain of discontinuity for $\Gamma_K$ and let
$\Sigma_K = \Omega_K/\Gamma_K$ be the quotient projective
structure. Since $\Omega_K \supseteq \Omega$, $D(\tilde{c}_K)$ will
also descend to two simple closed curves on $\Sigma_K$. We abuse
notation by also referring to these curves as $c_1$ and $c_2$.

\begin{lemma}
\label{surfschott}
The projective structure $\Sigma_K$ is homeomorphic to a surface of
genus $2g-1$. Furthermore there is an orientation reversing
involution
$\phi: \Sigma_K \longrightarrow \Sigma_K$ that fixes $c_1$ and $c_2$
pointwise and lifts to an orientation reversing, $\Gamma_K$-invariant
involution $\tilde{\phi}: \Omega_K \longrightarrow \Omega_K$ which
fixes $D(\tilde{c}_K)$ pointwise. 
\end{lemma}

{\bf Proof.} We reserve the 3-dimensional proof of this lemma to the
next section where we will prove the stronger Lemma \ref{schottky}. 
\qed{surfschott}

Since $D(\tilde{K})$ is contained in $\Omega_K$, $D(\tilde{K})/\Gamma$ is a subsurface of $\Sigma_K$.
Let $\Sigma^-_K$ be the closure of $D(K)/\Gamma_K$ in $\Sigma_K$. Then $\Sigma_K^-$ is homeomorphic to a genus $g-1$ surface with two boundary components $c_1$ and $c_2$. Let $\Sigma^+_K$ be the closure of the 
complement of $\Sigma^-_K$ in $\Sigma_K$. The involution $\phi$ from Lemma \ref{surfschott} will then restrict to a homeomorphism from $\Sigma^-_K$ to $\Sigma^+_K$ so $\Sigma^+_K$ is also a genus $g-1$ surface with two boundary components.

We also know that $D(\tilde{K})$ is contained in $\Omega$ so $\Sigma^-_K$ is also a subsurface of the cover $\Omega/\Gamma_K$ of $\Sigma$. In fact the covering map $\pi : \Omega/\Gamma_K \longrightarrow \Sigma$ restricts to a one-to-one map from the interior of $\Sigma^-_K$ to $\Sigma - c$ and is a two-to-one map from $c_1 \cup c_2$ to $c$.  We use $\pi$ to
define an equivalence relation for points $p_1 \in c_1$ and $p_2 \in
c_2$ with $p_1 \sim p_2$ if $\pi(p_1) = \pi(p_2)$. Then the quotient $\Sigma^-_K/\sim$ is exactly the original projective structure $\Sigma$. More importantly, the quotient $\Sigma_c = \Sigma^+_K/\sim$ will also be projective structure on $S$.

\begin{theorem}
\label{newproj}
$\Sigma_c$ is a projective structure on $S$ with holonomy $\rho$.
\end{theorem}

{\bf Proof.} We can explicitly write down a formula for a developing map for $\Sigma_c$ by modifying the developing map $D$ for $\Sigma$. Namely define $D_c: \tilde{S} \longrightarrow \chat$ by the formula
$$D_c(x) = \rho(g^{-1}) \circ \tilde{\phi} \circ D(g(x)) \mbox{\ if $g(x)$ is in the closure of $\tilde{K}$}.$$
It is a simple matter of retracing
definitions to see that $D_c$ is well defined, a developing map for
$\Sigma_c$ and has holonomy $\rho$. \qed{newproj}

\begin{cor}
\label{notgraft}
The projective structure $\Sigma_c$ is not obtained by grafting
$\Sigma$.
\end{cor}

{\bf Proof.} The developing map $D_c$ has the opposite orientation as
that of $D$ so $\Sigma_c$ cannot be a grafting of
$\Sigma$. \qed{notgraft} 

In the above work we have assumed that $c$ is non-separating. This is
not essential. In fact, after minor modifications, the construction
works for any collection $\cC$ of $n$ disjoint, homotopically distinct
and essential simple closed curves. If $\tilde{\cC}$ is the pre-image
of $\cC$ in $\tilde{S}$ then the action of $\pi_1(S)$ on $\tilde{S} - \tilde{\cC}$ will have $k$ orbits where $k$ is the number of components of $S - \cC$. We choose a component $\tilde{K}_i$ corresponding to each orbit and let $\Gamma_{K_i}$ be the subgroup of $\Gamma$ that fixes $D(\tilde{K}_i)$ with $\Omega_{K_i}$ the domain of discontinuity of $\Gamma_{K_i}$. Each projective structure $\Sigma_{K_i} = \Omega_{K_i}/\Gamma_{K_i}$ can then be cut into two pieces $\Sigma^-_{K_i}$ and $\Sigma^+_{K_i}$ and there is an involution $\phi_i$ of $\Sigma_{K_i}$ swapping the two pieces. Then the $\Sigma^-_{K_i}$ can be glued together to reform $\Sigma$. The $\Sigma^+_{K_i}$ can also be glued together to form a new projective structure $\Sigma_{\cC}$. As before we can explicitly define a developing map $D_{\cC} : \tilde{S} \longrightarrow \chat$ for $\Sigma_{\cC}$ by the formula
$$D_{\cC}(x) = \rho(g^{-1}) \circ \tilde{\phi}_i \circ D(g(x)) \mbox{\ if $g(x)$ is in the closure of $\tilde{K}_i$.}$$
Again, it is a simple matter of tracing through the definitions to see that $D_{\cC}$ is a developing map for a projective structure on $S$ and that the holonomy of $D_{\cC}$ is $\rho$.

We also remark that if $c$ is a component of $\cC$ and $\cC' = \cC - c$. Then
$\Sigma_{\cC}$ can also be obtained by either grafting $\Sigma_{\cC'}$
along $c$ or grafting $\Sigma_c$ along $\cC'$.

This construction also works if $\Gamma$ is quasi-fuchsian. In this
case we have two initial projective structures $\Sigma^-$ and
$\Sigma^+$ corresponding to the two components of the domain of
discontinuity. We leave the following theorem as an exercise for the
reader.

\begin{theorem}
The projective structure $\Sigma^-_{\cC}$ is equivalent to grafting
$\Sigma^+$ along $\cC$.
\end{theorem}

This leads us to make the following conjecture for projective
structures with singly degenerate holonomy:

\begin{conj}
Every projective structure with holonomy a singly degenerate group
$\Gamma$ is either:
\begin{enumerate}
\item $\Sigma$

\item $\Sigma_{\cC}$ for some collection $\cC$

\item grafting of $\Sigma$

\item grafting of $\Sigma_{\cC}$ along $\cC$
\end{enumerate}
\end{conj}

\section{Cone manifolds}

We carry over our notation from the previous section. Let $M =
(\hthree \cup \Omega)/\Gamma$ be the quotient 3-manifold with
boundary. By Bonahon's theorem (Theorem \ref{bonahon}), $M$ is
homeomorphic to $S \times [-1,1)$. 
The interior of $M$ will have a complete hyperbolic structure while
the boundary $S \times \{-1\}$ is the projective structure $\Sigma$.

We recall the construction described in the introduction, adding more
details. Let $c$ be an essential simple closed curve on $S$ and make
the further assumption that $c \times \{0\}$ is a geodesic in $M$. Let
$A = c \times [0,1)$ be an annulus in $M$. Then $A$ lifts homeomorphically to an annulus $A_\integers$ in the
$\integers$-cover $M_\integers$ of $M$ associated to $c$. Let
$\overline{M - A}$ and $\overline{M_\integers - A_\integers}$ be the
metric completions of $M - A$ and $M_\integers - A_\integers$,
respectively.

The boundaries of both $\overline{M - A}$ and $\overline{M_\integers -
  A_\integers}$ are isometric to two copies of $A$ glued at $c \times
\{0\}$. Orient $A$. We then distinguish between the two copies of $A$
in the boundary of $\overline{M - A}$ by labeling $A^+$ the copy of
$A$ where the normal points outward and $A^-$ the copy of $A$ where
the normal points inward. Similarly label the two copies of $A$ in the
boundary of $\overline{M_\integers - A_\integers}$, $A^+_\integers$
and $A^-_\integers$. All four of these annuli are isometric to $A$ and
we use this isometry to define an equivalence relation between points
on $A^+$ and $A^-_\integers$ and between $A^-$ and
$A^+_\integers$. Namely, if $p_1 \in A^+$ and $p_2 \in A^-_\integers$
then $p_1 \sim p_2$ if they are mapped to the same point in the
isometry to $A$. Similarly define an equivalence relation for points in
$A^-$ and $A^+_\integers$. Then
$$M_c = (\overline{M - A} \cup \overline{M_\integers - A_\integers})/\sim.$$

The hyperbolic structures on $M - A$ and on $M_\integers -
A_\integers$ will extend to a smooth hyperbolic structure in $M_c$
except at $c \times \{0\}$. At $c \times \{0\}$ the metric has a cone
singularity of cone angle $4\pi$. Furthermore $M_c$ is homeomorphic to
$S \times [-1,1)$ with $S \times \{-1\}$ the projective structure $\Sigma$. Our goal for the remainder of this section is to
show that $M_c$ is a quasi-fuchsian cone-manifold. That is we will show that $M_c$ extends to the projective structure $\Sigma_c$ on $S \times \{1\}$.

As in the previous section we assume for simplicity that $c$ is
non-separating. The general case is the same with more
notation. Let $B = c \times [-1,0]$ be an
annulus in $M$ and let $\tilde{B}_K = \tilde{c}_K \times [-1,0]$ be
the components of the pre-image of $B$ that bound $\tilde{K} \times [-1,0]$ in $\tilde{M}$. Let
$\tilde{L} = (\tilde{K} \times \{0\}) \cup \tilde{B}_K$.

Let $H = (\hthree \cup \Omega_K)/\Gamma_K$. Since $\Gamma_K$ restricts to an action on 
$\tilde{L}$, the quotient $L = \tilde{L}/\Gamma_K$ is a surface in $H$.
\begin{lemma}
\label{schottky}
$H$ is a genus $2g-1$ handlebody with boundary. Furthermore there is
an orientation reversing involution $\phi:H \longrightarrow H$
with $\phi|_L \equiv \id$ which lifts to an orientation reversing,
$\Gamma_K$-equivariant involution $\tilde{\phi}: (\hthree \cup
\Omega_K) \longrightarrow (\hthree \cup \Omega_K)$ with
$\tilde{\phi}|_{\tilde{L}} \equiv \id$.
\end{lemma}

{\bf Proof.} The covering map $H \longrightarrow M$ is infinite-to-one. In particular, on the single end of $H$ it is infinite-to-one. 
By the covering theorem (\cite{Canary:inj:radius}),
either $\Gamma_K$ is geometrically 
finite or $M$ has a finite index cover that fibers over the circle. Since the latter condiation is not true $\Gamma_K$ must be geometrically finite. Furthermore, $\Gamma_K$ does not contain parabolics. Therefore $\Gamma_K$ is
a Schottky group with $2g-1$ generators and $H = (\hthree \cup
\Omega_K)/\Gamma_K$ is a genus $2g-1$ handlebody with boundary.

The inclusion of $L$ in $H$ is a homotopy equivalence. Therefore $H$
is homeomorphic to $S' \times [-1,1]$ 
where $S'$ is a genus $g-1$ surface with two boundary components and
$S' \times \{0\} = L$. This product structure defines an obvious
involution of $H$ which lifts to the universal cover to obtain the
desired involution $\tilde{\phi}$ of $\hthree \cup \Omega_K$. \qed{schottky}

{\bf Remark.} Note that although the handlebody $H$ covers $M$ the
product structure we have chosen for $H$ is not equivariant and does
not descend to the product structure on $M$.

\medskip

\begin{theorem}
\label{coneman}
The hyperbolic cone-manifold $M_c$ is quasi-fuchsian with projective
boundary $\Sigma$ and $\Sigma_c$.
\end{theorem}

{\bf Proof.} To prove the theorem we make an alternative construction of $M_c$.

We begin with on observation about the surface $S$. Let $S'$ be the cover of $S$ corresponding to $\pi_1(S - c)$. As we have already noted $c$ has two homeomorphic lifts $c_1$ and $c_2$. Next we divide $S'$ into three subsurfaces $S_0$, $S_1$ and $S_2$ with $S_0$ a compact genus $g-1$ surface with two boundary components and $S_1$ and $S_2$ both homeomorphic to the annulus $S^1 \times [0,1)$. We also assume that $S_0 \cap S_1 = c_1$ and $S_0 \cap S_2 = c_2$. Note that the covering map $\pi: S' \longrightarrow S$ defines an equivalence relation on points $p_1 \in c_1$ and $p_2 \in c_2$ by $p_1 \sim p_2$ if $\pi(p_1) = \pi(p_2)$. Then $\pi$ restricts to a homeomorphism from the quotient $S_0/\sim$ to $S$. On the quotient $(S_1 \cup S_2)/\sim$, $\pi$ becomes the covering map for the $\integers$-cover of $S$ associated to $c$.

Next take the product $M = S \times (-1,1)$ and let $X_i = S_i \times (-1,1)$. Extending our equivalence relation to the product structure in the obvious way we then see that $X_0/\sim$ is homeomorphic and isometric to $M$ while $(X_1 \cup X_2)/\sim$ is the $\integers$-cover $M_\integers$ of $M$ associated to $c \times \{0\}$.

\begin{figure}[h]
\begin{center}
\epsfig{file=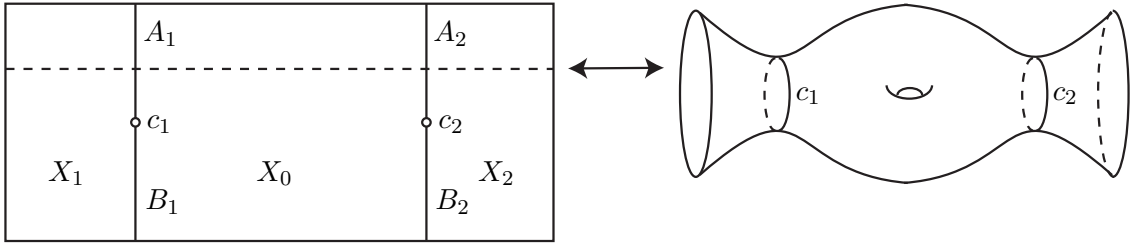, scale=.85}
\caption{The rectangle gives a schematic picture of the product structure on $H$. The horizontal lines represent the cover $S'$ of $S$}
\label{rectangle}
\end{center}
\end{figure}

To construct $M_c$ we subdivide the annuli that bound the $X_i$. The boundary of $X_1$ is the annulus $c_1 \times (-1,1)$. Let $A^+_1 = c_1 \times [0,1)$ and $B^+_1 = c_1 \times (-1,0]$. Similarly divide the boundary of $X_2$ into two annuli $A^-_2$ and $B^-_2$. We also divide each of the two annuli that bound $X_0$ into two sub-annuli $A^-_1$, $B^-_1$, $A^+_2$ and $B^+_2$. To construct $\overline{M-A}$ we start with $X_0$ and glue $B^-_1$ to $B^+_2$. To construct $\overline{M_\integers - A_\integers}$ we glue $X_1$ to $X_2$ by attaching $B^+_1$ to $B^-_2$. Finally, to construct $M_c$ we glue the $A$ annuli together. Namely we glue $A^+_1$ to $A^-_1$ and $A^+_2$ to $A^-_2$.

Of course this is simply restating our original construction of $M_c$. As an alternative we first glue the $A$ annuli and then glue the $B$ annuli. In both cases we use the same gluing pattern so we  get the same hyperbolic structure $M_c$. To see the advantage of gluing in this order we recall that the cover $S' \times (-1,1)$ of $M$ is the interior of the handlebody $H$. The boundary of $H$ is the  projective structure $\Sigma_K$. The annulus $B$ lifts to two annuli $B_1$ and $B_2$ in $H$ which extend to closed curves $c_1$ and $c_2$ on $\Sigma_K$. Next we note that when we glue $X_1$ and $X_2$ to $X_0$ along the $A$ annuli we get the metric completion of $H - (B_1 \cup B_2)$. This compact manifold has boundary consisting of the $B$ annuli and the projective structures $\Sigma^+_K$ and $\Sigma^-_K$. When we glue the $B$ annuli the two boundary curves of $\Sigma^+_K$ are identified to form the projective structure $\Sigma_c$. Similarly the boundary curves of $\Sigma^-_K$ are identified to form the original projective structure $\Sigma$. Therefore $M_c$ is compactified by its projective boundary and is a quasi-fuchsian cone-manifold.
 \qed{coneman}

\begin{figure}[h]
\begin{center}
\epsfig{file=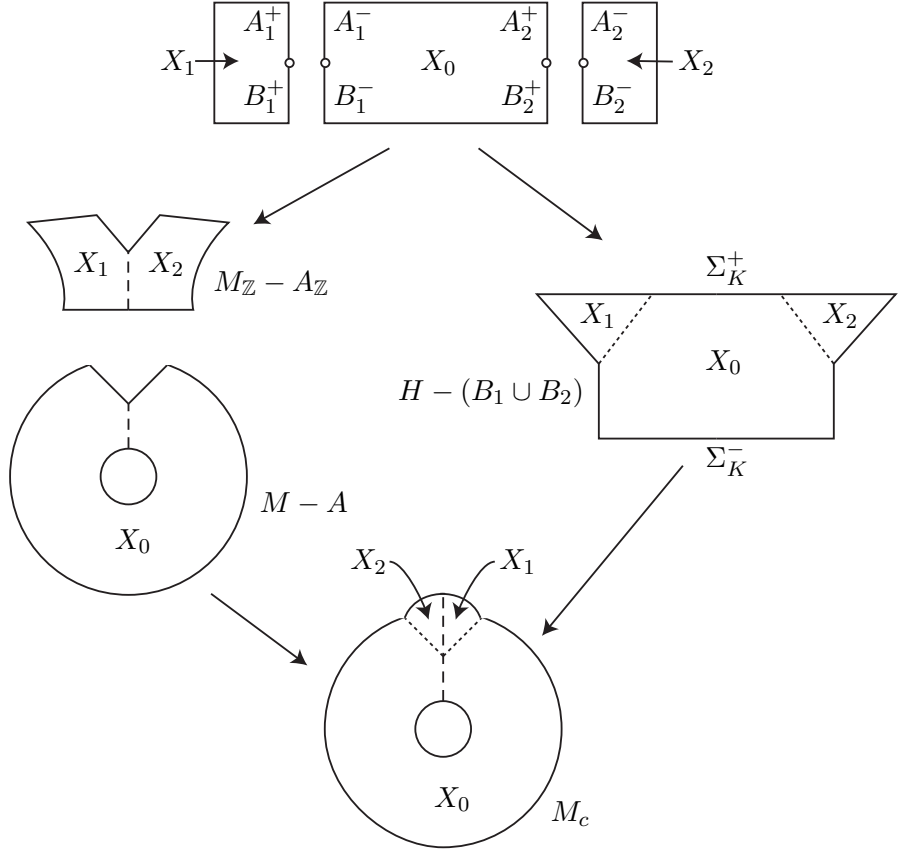, scale=.7}
\caption{The figure gives a schematic description of the two constructions of $M_c$. On the left is the original construction while on the right is the alternative construction.}
\end{center}
\end{figure}

\section{The Bers' conjecture}
\label{unboundedbers}
In the previous section we constructed quasi-fuchsian hyperbolic
cone-manifolds. We now use the deformation theory of hyperbolic
cone-manifolds to show that these cone structures are geometrically
close to a smooth quasi-fuchsian structure. The analytic deformation theory of hyperbolic cone-manifolds was developed by Hodgson and Kerckhoff in a series of papers (\cite{Hodgson:Kerckhoff:cone, Hodgson:Kerckhoff:dehn, Hodgson:Kerckhoff:tube}) and extended to the geometrically finite setting in \cite{Bromberg:rigidity, Bromberg:schwarz}. The basic idea is that if
the cone singularity is short and has a large tube radius then there
is a one parameter family of cone-manifolds decreasing the cone angle
from $4\pi$ to a cone-manifold with cone angle $2\pi$. When the cone
angle is $2\pi$ the hyperbolic structure is non-singular.

Although the theory applies in greater generality, we will confine
ourselves to a quasi-fuchsian cone-manifolds. The following result is essentially Theorems 1.2 and 1.3 of \cite{Bromberg:schwarz}.

\begin{theorem}
\label{schwarzbound}
Suppose $M_\alpha$ is a quasi-fuchsian cone manifold with cone singularity $c$, cone angle $\alpha$ and conformal boundary $X$ and $Y$. Also assume the tube radius of $c$ is greater than $\sinh^{-1}\sqrt{2}$. Then:
\begin{enumerate}
\item There exist an $\ell_0>0$ depending only on $\alpha$ such that for all $t \leq \alpha$ there exists a quasi-fuchsian cone-manifold $M_t$ with cone singularity $c$, cone angle $t$ and conformal boundary $X$ and $Y$.

\item Furthermore if $\Sigma_\alpha$ and $\Sigma_t$ are the projective boundaries corresponding to $X$ for $M_\alpha$ and $M_t$, respectively, there exists a $K$ depending only on $\alpha$, $\|\Sigma_\alpha\|_\infty$ and the injectivity radius of the hyperbolic metric on $X$ such that
$$d(\Sigma_\alpha, \Sigma_t) \leq K \operatorname{length}(c)$$
where the length is measured in the $M_\alpha$-metric.
\end{enumerate}
\end{theorem}

We can now prove our main theorem:

\begin{theorem}
\label{unbounded}
Assume that $\Gamma \in AH(S)$ has no parabolics.
If $M = \hthree/\Gamma$ is singly degenerate and has unbounded
geometry then we have $\Gamma \in \overline{B}_X$ where $X$ is the conformal
boundary of $M$.
\end{theorem}

{\bf Proof.} By the Margulis lemma there exists an $\ell_1$ such that
if $c$ is closed geodesic in $M$ with $\operatorname{length}(c)<
\ell_1$ then $c$ has an
embedded tubular neighborhood of radius
$\sinh^{-1}\sqrt{2}$. We need the 
following theorem of Otal:
\begin{theorem}[Otal \cite{Otal:short}]
\label{otal}
Let $c$ be a simple closed geodesic in $M$. There exists an $\ell_2
>0$ such that if $\operatorname{length}(c) < \ell_2$ then $c$ is isotopic to a
simple closed curve on $S \times \{0\}$ in $M$.
\end{theorem}
Let $\ell = \min(\ell_0, \ell_1, \ell_2)$ where $\ell_0$ is the
constant from Theorem \ref{schwarzbound}.

Since $M$ has unbounded geometry there are a sequence of closed
geodesics $c_i$ in $M$ with $\operatorname{length}(c_i) \rightarrow
0$. Therefore we can assume that $\operatorname{length}(c_i) < \ell$ for all 
$i$. We can then apply Theorem \ref{coneman} to construct a sequence
of cone-manifolds $M_i$ with cone-singularity $c_i$ and cone-angle
$4\pi$. Furthermore, an embedded tubular neighborhood of $c_i$ in
$M$ will lift to an embedded tubular neighborhood of $c_i$ in $M_i$ of
the same radius. Therefore $c_i$ will have an embedded tubular
neighborhood of radius $\sinh^{-1}\sqrt{2}$ in
$M_i$. 

We can now apply Theorem \ref{schwarzbound} to the $M_i$. 
If $X$ and $Y_i$ are the components of conformal boundary of $M_i$ let
$M'_i$ be the quasi-fuchsian cone manifold with cone singularity $c_i$, cone angle $2\pi$ and
conformal boundary $X$ and $Y_i$ given by (a) of Theorem \ref{schwarzbound}. Since the cone angle is $2\pi$ the hyperbolic structure on $M'_i$ will be smooth so there will be a unique Kleinian group $\Gamma_i \in B_X$ such that $M'_i = \hthree/\Gamma_i$. Note that for each $M_i$ the
component of the projective boundary associated to $X$ will be
$\Sigma$, the projective boundary of the original hyperbolic
structure $M$. Let $\Sigma_i$ be the component of the projective
boundary of $M'_i$ associated to $X$. By Theorem $\ref{schwarzbound}$, 
$$d(\Sigma, \Sigma_i) \leq K\operatorname{length}(c).$$
Therefore we have $\Sigma_i \rightarrow \Sigma$ in $P(X)$ which implies that $\Gamma_i
\rightarrow \Gamma$ in $AH(S)$. Since each $\Gamma_i$ is contained in
$B_X$, we conclude $\Gamma \in \overline{B}_X$. \qed{bers}

Combining Theorem \ref{unbounded} with Theorem \ref{minsky} we have:

\begin{theorem}
\label{bers}
Assume that $\Gamma \in AH(S)$ has no parabolics. If $M = \hthree/\Gamma$ is singly degenerate then $\Gamma \in \overline{B}_X$ where $X$ is the conformal boundary of $M$.
\end{theorem}

\bibliographystyle{../tex/math}
\bibliography{../tex/math}

\begin{sc}
\noindent
Department of Mathematics\\
California Institute of Technology\\
Mailcode 253-37\\
Pasadena, CA 91125\\
\end{sc}

\end{document}